\documentclass{article}

     \usepackage[preprint]{neurips_2020}

\usepackage[utf8]{inputenc} 
\usepackage[T1]{fontenc}    
\usepackage{hyperref}       
\usepackage{url}            
\usepackage{booktabs}       
\usepackage{amsfonts}       
\usepackage{nicefrac}       
\usepackage{microtype}      

\usepackage{hyperref}
    \hypersetup{colorlinks=true, linkcolor=black, citecolor=black, urlcolor=black}

\usepackage{microtype}

\usepackage[shortlabels, inline]{enumitem} 

\usepackage[retainorgcmds]{IEEEtrantools}

\usepackage{algorithm}
\usepackage{algorithmic}

\usepackage[normalem]{ulem}

\usepackage{natbib}

\usepackage{amssymb}
\usepackage{amsmath}
\usepackage{amsthm}
\usepackage{bbm} 
	\newcommand{\1}{\mathbbm{1}}
    \providecommand{\1}{\mathds{1}}

\usepackage{pgfplots}
\usepackage{pgfplotstable}
\usepackage{tikz}
    \definecolor{myblue}{HTML}{3333CC}
    \definecolor{mylightblue}{HTML}{CCCCFF}
    \definecolor{myorange}{HTML}{FF6600}
    \definecolor{myred}{HTML}{CC3333}
    \definecolor{mygreen}{HTML}{11AA11}
    \definecolor{mylightgray}{HTML}{E0E0E0}
    \usetikzlibrary{decorations.pathmorphing}
    \tikzstyle{vertex}=[circle, draw, fill=white, inner sep=0pt, minimum width=1ex]
    \tikzset{every picture/.append style={baseline,scale=1.1}}
    \tikzset{every edge/.append style={thick}}     
    \tikzstyle{cut-edge}=[dotted]
    \tikzstyle{base}=[line width=1.2pt]
    \tikzstyle{lifted}=[blue]

\usetikzlibrary{external}

\providecommand{\R}{\mathbb{R}}
\providecommand{\Z}{\mathbb{Z}}
\providecommand{\N}{\mathbb{N}}

\providecommand{\B}{\mathcal{B}}
\providecommand{\C}{\mathcal{C}}
\providecommand{\D}{\mathcal{D}}

\providecommand{\F}{\mathcal{F}}

\renewcommand{\P}{\mathcal{P}}
\providecommand{\M}{\mathcal{M}}

\providecommand{\abs}[1]{\lvert #1 \rvert}

\theoremstyle{plain}
\newtheorem{thm}{Theorem}
\newtheorem{lem}{Lemma}

\newtheorem{cor}{Corollary}

\theoremstyle{definition}
\newtheorem{defn}{Definition}

\newtheorem{ex}{Example}

\theoremstyle{remark}

\DeclareMathOperator{\conv}{conv}

\providecommand{\mc}{\mathsf{MC}}
\providecommand{\cut}{\mathsf{CUT}}

\providecommand{\pflow}{\mathsf{P}_{\mathrm{flow}}}

\title{Flow-Partitionable Signed Graphs}

\author{%
  Jan-Hendrik Lange \\
  Max Planck Institute for Informatics \\
  Saarbr\"ucken, Germany \\
}

\begin{document}

\maketitle

\begin{abstract}
The NP-hard problem of correlation clustering is to partition a signed graph such that the number of conflicts between the partition and the signature of the graph is minimized.
This paper studies graph signatures that allow the optimal partition to be found efficiently.
We define the class of flow-partitionable signed graphs, which have the property that the standard linear programming relaxation based on so-called cycle inequalities is tight.
In other words, flow-partitionable signed graphs satisfy an exact max-multiflow-min-multicut relation in the associated instances of minimum multicut.
In this work we propose to characterize flow-partitionable signed graphs in terms of forbidden minors.
Our initial results include two infinite classes of forbidden minors, which are sufficient if the positive subgraph is a circuit or a tree.
For the general case we present another forbidden minor and point out a connection to open problems in the theory of ideal clutters.
\end{abstract}


\section{Introduction}

In light of the good performance of linear programming (LP) relaxations for NP-hard combinatorial optimization problems in structured prediction, there is considerable research interest in theoretical explanations for this phenomenon.
A natural question in this context is the following: For which inputs does the solution of the LP relaxation coincide with the solution of the unrelaxed problem?
In this paper we propose an analysis of such a tightness property for the standard LP relaxation of the correlation clustering problem.

Correlation clustering is the problem of finding a partition of a signed graph that minimizes the number of \emph{errors}, which arise when either a negative edge is within a cluster or a positive edge is between clusters \citep{Bansal2004}.
It is encountered under different names in other communities due to slightly altered formulations, such as \emph{minimum multicut} problem.
When only partitions into two components are considered, the problem is called 2-correlation clustering and is closely related to the max-cut problem and thus binary quadratic programming.
The two problems differ from other clustering formulations by their reliance on qualitative edge information, which can be interpreted as attraction or repulsion between entities in a network depending on the edge sign.
As fundamental problems in signed graph partitioning they have found applications in a wide variety of areas such as probabilistic graphical models \citep{Wainwright2008}, statistical physics \citep{Liers2005}, social network analysis \citep{Cesa-Bianchi2012,Veldt2018}, image analysis \citep{Andres2011,Kappes2011,Keuper2015b,Beier2017} and computer vision \citep{Insafutdinov2016,Keuper2015a} to name a few.

Despite a lot of progress in approximation algorithms for correlation clustering \citep{Chawla2015}, the sign patterns that prohibit \emph{optimal} solutions to be found efficiently are not well understood.
Signed graphs that can be clustered \emph{without errors} into two components or an arbitrary number of components are called balanced \citep{Harary1953}, respectively weakly balanced \citep{Davis1967}, and are characterized by the absence of certain signed cycles.
Moreover, in the case of 2-correlation clustering, a precise characterization of the tightness of a standard LP relaxation in terms of forbidden minors is available \citep{Guenin2001,Weller2016a}.
In this paper we propose to extend this analysis to correlation clustering. 
Namely, we study the structure of signed graphs that allow an optimal solution of the associated correlation clustering problem to be found efficiently by means of linear programming.
More precisely, we introduce the class of \emph{flow-partitionable} signed graphs, for which the standard relaxation based on so-called cycle inequalities is tight, see Section \ref{sec:signed-graphs}.
By reduction to minimum multicut, the class of flow-partitionable signed graphs define instances for which an exact max-multiflow-min-multicut relation holds.
Furthermore, we propose to characterize the class of flow-partitionable signed graphs in terms of forbidden minors via the correspondence between signed graphs and the so-called flow clutter, as described in Section~\ref{sec:ideal-clutters}.
We present two infinite classes of forbidden minors which provide a complete characterization of tightness for the special cases that the positive subgraph is either a circuit or a tree.
For the general case we describe another forbidden minor and point out a connection to open problems in the theory of ideal clutters, which suggests that a complete characterization is more challenging, see Section~\ref{sec:idealness}.
Proofs for our results can be found in the Appendix.
Our work enhances the understanding of the performance of LP relaxations in correlation clustering and thus, more generally, in signed graph partitioning and structured prediction.


\section{Related Work}
\label{sec:related-work}

Signed graphs have been introduced in the study of social networks to identify coherent social groups \citep{Cartwright1956}.
\cite{Harary1953} and \cite{Davis1967} showed that a signed graph can be partitioned into exactly two or arbitrarily many components without errors if it does not contain any cycles with an odd number of negative edges, respectively exactly one negative edge.
The signed graphs that satisfy the former or the latter property are called \emph{balanced}, respectively \emph{weakly balanced}.

The problem of partitioning into two components with as few errors as possible is well-studied from a polyhedral perspective \citep{Groetschel1981,Barahona1983,Barahona1986}.
\cite{Groetschel1981} introduced the class of \emph{weakly bipartite} signed graphs, which have the property that the linear programming relaxation defined by odd cycle inequalities is tight.
\cite{Guenin2001} showed that a signed graph is weakly bipartite if, and only if, it has no odd $K_5$ signed minor, which is a seminal result in combinatorial optimization and generalizes earlier findings.
The characterization can be translated into a similar tightness condition for binary quadratic programming problems and thus binary graphical models \citep{Weller2016a,Michini2016}.

For the problem of partitioning optimally into an arbitrary number of components, polyhedral works include \citep{Groetschel1989,Groetschel1990,Deza1990, Deza1992,Chopra1993}.
Further, \cite{Chopra1994} shows that if the graph has treewidth at most two, then the multicut polytope is fully described by cycle inequalities, which means any signed graph without a $K_4$ minor is flow-partitionable.
In the machine learning community the problem is known as correlation clustering and has been mainly studied from an approximation perspective.
The name correlation clustering is due to the interpretation that the optimal clustering maximally correlates with the graph signature.
Hardness results and approximation algorithms for particular classes of graphs and/or edge weights are due to
\cite{Bansal2004, Charikar2005, Demaine2006, Chawla2006, Chawla2015, Ailon2008, Ailon2012, Klein2015, Veldt2017}.
Heuristic methods that determine a clustering by greedily contracting edges are proposed by \cite{Keuper2015b,Levinkov2017, Kardoost2019,Bailoni2019}.
Methods that solve a Lagrangian relaxation of the problem are due to \cite{Yarkony2012,Yarkony2015,Swoboda2017}.
\cite{Lange2019} develop combinatorial criteria that allow to identify parts of optimal solutions efficiently.

The correlation clustering problem can be reduced to minimum multicut by substituting every negative edge with a positive edge and a terminal pair \citep{Demaine2006}.
Thus, a flow-partitionable signed graph defines an instance of the minimum multicut problem that satisfies an exact max-multiflow-min-multicut relation.
Although there is a large body of work on multi-commodity flows and multicuts \citep{Schrijver2003}, commonly the focus has been on conditions that guarantee the existence of an optimal integer multi-commodity flow.
\cite{Cornaz2011} raised the question when an exact min-max-relation between the minimum multicut and integer maximum multi-commodity flow holds.
Furthermore, he conjectures a characterization in terms of forbidden strong minors.
In contrast, we consider the case when the path relaxation of the multicut problem gives an exact solution, regardless of how the dual problem behaves.
We adopt some of the terminology of \cite{Cornaz2011} for our approach and characterize flow-partitionable signed graphs for the cases when the positive subgraph is either a tree or a circuit.
While the former case corresponds to multicut in trees, which is NP-hard \citep{Garg1997}, in the latter case correlation clustering can be solved in polynomial time, similar to multicut \citep{Bentz2009}. 
The path relaxation for multicut is integral for up to two terminal pairs, which is a consequence of the well-known \emph{max-flow-min-cut}-Theorem \citep{Ford1956}, respectively a result due to \cite{Hu1963} on two-commodity flows, cf. \cite[Cor71.1d]{Schrijver2003}.
These results carry over to correlation clustering and imply that signed graphs are flow-partitionable for up to two negative edges.
Other special cases that do not directly carry over are due to \cite{Tang1964, Sakarovitch1966, Rothfarb1969, Karzanov1989}.
To the best of our knowledge, a general characterization of flow-partionable signed graphs is an open problem.

Our work relies on results in the theory of ideal clutters \citep{Lehman1979,Lehman1990,Seymour1990, Padberg1993,Cornuejols1994}, in particular Lehman’s theorem, which is a key ingredient towards proving integrality of LP relaxations of covering problems.
Other applications of Lehman's theorem include \citep{Guenin2001, Schrijver2002, Abdi2016}.


\section{Correlation Clustering and Flow-Partitionable Signed Graphs}
\label{sec:signed-graphs}

In this paper we consider signed graphs $G = (V,E)$, where each edge $e \in E$ is labeled either~$+$ or~$-$.
We collect the respective edges in two sets $E^+$ and $E^-$, which we call the \emph{positive}, respectively \emph{negative} edges and thus $E = E^+ \cup E^-$.
In order to emphasize that $G$ is signed, we may write $G = (V,E^+,E^-)$ instead of $G= (V,E)$.
Further, let $G^+ = (V, E^+)$ and $G^- = (V, E^-)$ denote the positive, respectively negative, subgraph of $G$.

A \emph{graph partition} or \emph{clustering} of $G$ is a set $\Pi = \{ U_1, \dotsc, U_k \}$, where $k$ is arbitrary, such that $U_i \subseteq V$, $V = \bigcup_i U_i$ and $U_i \cap U_j = \emptyset$ for all $i \neq j$.
The $U_i$ are called \emph{components} of the partition or \emph{clusters}.
Every graph partition is associated with the set of edges within, respectively between the components of the partition. The latter set is also known as \emph{multicut}.

\begin{defn}[Multicut/Cut]
Let $\Pi = \{ U_1, \dotsc, U_k \}$ be any partition of the graph $G = (V,E)$.
The set of edges between components
\begin{align*}
M = \bigcup_{1 \leq i < j \leq k} \delta(U_i, U_j)
\end{align*}
is called \emph{multicut} of $G$ associated with $\Pi$. 
If $k=2$, then $M = \delta(U_1,U_2) = \delta(U_1) = \delta(U_2)$ is also called a \emph{cut} of $G$.
\end{defn}

\begin{defn}[Multicut/Cut Polytope]
The convex hull of characteristic vectors of multicuts, respectively cuts of $G$ is called \emph{multicut polytope}, respectively \emph{cut polytope} of $G$ and is denoted by
\begin{align*}
\mc(G) = \conv \big \{ \1_M \mid M \text{ multicut of } G \big \}, \qquad \cut(G) = \conv \{ \1_{\delta(U)} \mid U \subseteq V \}.
\end{align*}
We write $\mc = \mc(G)$, respectively $\cut = \cut(G)$ for short.
\end{defn}

For any clustering of a signed graph $G = (V, E^+, E^-)$, we say that positive edges between clusters and negative edges within clusters constitute \emph{errors} of the clustering w.r.t.\ the signature of the graph.
\begin{defn}[Correlation Clustering]
The \emph{correlation clustering} problem is to find a clustering of signed graph $G=(V,E^+,E^-)$ that minimizes the number of errors and can be formulated as
\begin{align}
\min_{x \in \mc} \; \sum_{e \in E^-} (1-x_e) + \sum_{e \in E^+} x_e. \label{eq:correlation-clustering} \tag{CC}
\end{align}
Similarly, if only partitions into two components are allowed, the problem is called \emph{2-correlation clustering} and is analogously formulated as
\begin{align}
\min_{x \in \cut} \; \sum_{e \in E^-} (1-x_e) + \sum_{e \in E^+} x_e. \label{eq:2-correlation-clustering} \tag{2CC}
\end{align}
\end{defn}

Problem \eqref{eq:correlation-clustering} can be reduced to the \emph{minimum multicut} problem by substitution of every negative with a positive edge and a terminal pair \citep{Demaine2006}.
Problem \eqref{eq:2-correlation-clustering} specializes to the \emph{max-cut} problem by setting $E^+ = \emptyset$.
Both problems \eqref{eq:correlation-clustering} and \eqref{eq:2-correlation-clustering} can be solved exactly by means of integer linear programming (ILP), which employs a linear relaxation of $\mc$, respectively $\cut$.
In its most basic form it includes inequalities associated with cycles of the graph \citep{Barahona1986,Chopra1993}.

It is straightforward to show, however, that it suffices to consider only a subset of all cycles, which depend on the signature of $G$.
For the purpose of this paper, a cycle of $G$ is a subset of edges $C \subset E$ such that any pair of (cyclically) adjacent edges shares a vertex.
A cycle without repeated vertices is called \emph{circuit}.
A circuit is called \emph{odd} if it contains an odd number of negative edges.
A circuit is called \emph{bad} if it contains exactly one negative edge.
\emph{Balanced} signed graphs have no odd circuits and \emph{weakly balanced} signed graphs have no bad circuits \citep{Harary1953, Davis1967}.
Naturally, any balanced signed graph is also weakly balanced.
For simplicity, and in line with \citep{Cornaz2011}, we shall also refer to a bad circuit as \emph{flow}.
Let $\F = \F(G)$ denote the set flows of the signed graph $G$.

For balanced or weakly balanced signed graphs, problem \eqref{eq:2-correlation-clustering}, respectively \eqref{eq:correlation-clustering}, is trivial.
Otherwise, an ILP formulation of \eqref{eq:correlation-clustering} (and analogously for \eqref{eq:2-correlation-clustering}) can be derived via the substitution $\hat x_e = 1 - x_e$ for $e \in E^-$ and $\hat x_e = x_e$ for $e \in E^+$ (see e.g.~\citep{Lange2018}), which yields
\begin{align}
\min \; \sum_{e \in E} \hat x_e \quad
\text{s.t.} \quad \sum_{e \in C} \hat x_e \geq 1 \; \forall C \in \F, \; \hat x \geq 0, \; \hat x \in \Z^E. \label{eq:flow-covering-problem}
\end{align}
An LP relaxation of \eqref{eq:flow-covering-problem} called \emph{cycle relaxation} is obtained by dropping the integrality constraints.
Similarly, we obtain an ILP formulation and a corresponding cycle relaxation of \eqref{eq:2-correlation-clustering} by substitution of $\F$ with the set of odd circuits.

The main motivation of this paper is the characterization of those signed graphs for which the cycle relaxation of \eqref{eq:flow-covering-problem} is tight, i.e.\ when the \emph{flow covering polyhedron}
\begin{align}
\pflow(G) = \left \{ x \geq 0 \; \Bigg | \sum_{e \in C} x_e \geq 1 \quad \forall C \in \F \right \} \label{eq:flow-covering-polyhedron}
\end{align}
is integral.
We shall refer to such signed graphs as \emph{flow-partitionable}.
\begin{defn}
A signed graph $G$ is called \emph{flow-partitionable} if its associated flow covering polyhedron $\pflow(G)$ is integral.
\end{defn}
Analogously, the signed graphs for which the \emph{odd circuit covering} polyhedron is integral are the  \emph{weakly bipartite} signed graphs \citep{Groetschel1981}.
They were characterized by \cite{Guenin2001} in terms of forbidden minors.
Our work is a first step towards the (arguably more challenging) characterization of flow-partitionable signed graphs in terms of forbidden minors.


Note that the integrality properties described here imply that the weighted versions of problems \eqref{eq:correlation-clustering} and \eqref{eq:2-correlation-clustering} can be solved via their cycle relaxations for any choice of weights given that the signature of the graph is fixed.


\section{Ideal Clutters and the Flow Clutter}
\label{sec:ideal-clutters}

In this section we frame the integrality of $\pflow(G)$ in the more general theory of ideal clutters.
Clutters are abstract objects that define covering polyhedra and idealness of the former corresponds to integrality of the latter.
Here, we introduce the necessary terminology from the theory of ideal clutters and define the flow clutter.

With any 0--1-matrix $A$ we associate the (fractional) covering polyhedron
\begin{align}
\mathsf P_A = \{ x \geq 0 \mid Ax \geq \1 \}.
\end{align}
The 0--1-matrix $A$ is called \emph{ideal} if $\mathsf P_A$ is integral.

\subsection{Clutters and Lehman's Theorem}

A family of subsets $\C$ of a finite ground set $E(\C)$ is called \emph{clutter} if no member of $\C$ is contained in another.
With any clutter $\C$ we naturally identify the 0--1-matrix $ A(\C)$ whose rows are the characteristic vectors of the members of $\C$.
By definition, no row vector of~$A(\C)$ dominates another, thus the constraint system that defines $\mathsf P_{A(\C)}$ is irredundant.
A clutter $\C$ is called \emph{ideal} if the associated matrix $A(\C)$ is ideal and otherwise it is called \emph{non-ideal}.

A \emph{cover} of $\C$ is a subset $B \subseteq E(\C)$ such that for all $C \in \C$ it holds that $B \cap C \neq \emptyset$.
The \emph{blocker} of $\C$ is the clutter that is formed by the minimal covers of $\C$.


The \emph{contraction} of an element $e \in E(\C)$ gives a clutter $\C / e$ over the ground set $E(\C) \setminus \{e\}$ that consists of the minimal sets from $\{ C \setminus \{e\} \mid C \in \C \}$.
The \emph{deletion} of an element $e \in E(\C)$ gives the clutter $\C \backslash e = \{ C \in \C \mid e \notin C \}$ over the ground set $E(\C) \setminus \{e\}$.
A \emph{minor} of $\C$ is any clutter obtained from $\C$ by a series of contraction and deletion operations.
Contraction and deletion operations are commutative.
The minor operations on clutters correspond naturally to restricting the polyhedron $\mathsf P_{A(\C)}$  by setting variables to $0$ or $1$.
More precisely, it holds that
\begin{align*}
P_{A(\C / e)} = P_{A(\C)} \cap \{ x \mid x_e = 0\}, \\
P_{A(\C \backslash e)} = P_{A(\C)} \cap \{ x \mid x_e = 1\}.
\end{align*}
The property of idealness is closed under taking minors \citep{Seymour1977}.

%

A clutter is called \emph{minimally non-ideal (MNI)} if it is not ideal but any (proper) minor is ideal.
Clearly, a clutter is ideal if, and only if, it has no MNI\ minor.

\begin{ex}
For any $2 \leq k \in \N$, the clutter $\D = \big \{ \{1, \dotsc, k\}, \{0,1\}, \dotsc, \{0,k\} \big \}$ over the ground set $\{0,1, \dotsc, k\}$ is called \emph{degenerate projective plane} of order $k$.
Any degenerate projective plane is MNI\ \citep{Lehman1979}.
\end{ex}

\paragraph{Lehman's Theorem}

For any MNI\ clutter $\C$, the \emph{core} of $\C$, denoted by $\overline \C$, is the clutter of minimum size members of $\C$.
\cite{Lehman1990} proved an important characterization of MNI clutters in terms of their cores.
We state Lehman's theorem below for later reference, an accessible proof can also be found in \citep{Seymour1990}. 

\begin{thm}[\cite{Lehman1990}] \label{thm:lehman}
Let $\C$ be an MNI\ clutter that is not a degenerate projective plane and let $\B$ be its blocker.
Then both $\overline \C = \{C_i\}_i$ and $\overline \B = \{B_i\}_i$ consist of $n = E(\C)$ members and can be ordered suitably such that for some $c,b \in \N$ it holds that
\begin{enumerate}[(i)]
\item $cb \geq n + 1$
\item $\forall C \in \overline \C : \abs{C} = c$ and $\forall B \in \overline \B : \abs{B} = b$
\item $\forall e \in E(\C) : \abs{\{C \in \overline \C \mid e \in C\}} = c$ and $ \abs{\{B \in \overline \B \mid e \in B\}} = b$
\item $\forall i,j \in [n] : \abs{C_i \cap B_j} = \begin{cases} cb - n + 1 & \text{if } i = j \\ 1 & \text{else} \end{cases}$
\item $\forall e,f \in E(\C) : \abs{\{i \in [n] : e \in C_i, f \in B_i\}} = \begin{cases} cb - n + 1 & \text{if } e = f \\ 1 & \text{else.} \end{cases}$
\end{enumerate}
In particular, the polyhedron $\mathsf P_{A(\C)}$ has the unique fractional vertex $\frac{1}{c} \1$.
\end{thm}

\subsection{The Flow Clutter}

Apparently, the set of flows $\F$ of a signed graph $G = (V, E^+, E^-)$, identified by their edge sets, is a clutter, which we call the \emph{flow clutter}.
We define a \emph{strong minor} \citep{Cornaz2011} of $\F$ as a minor (without singleton elements) that is obtained by a series of contractions of positive edges and deletions of arbitrary edges.
A \emph{strong minor} of the signed graph $G$ is a minor (without self-loops) that is obtained, analogously, by contraction of positive edges only and deletion of arbitrary edges.
Clearly, the strong minors are those minors that correspond to flow clutters defined by signed graphs.

\begin{cor}
Let $G = (V, E^+,E^-)$ be a signed graph and $\F$ its flow clutter.
Then any strong minor of $\F$ corresponds to a strong minor of $G$ and vice versa.
\end{cor}

We further call a flow clutter \emph{weakly MNI} if it is not ideal, but any (proper) strong minor is ideal.
It is clear that any MNI clutter is also weakly MNI.
However, weakly MNI clutters may have MNI minors that are constructed by contraction of negative edges.
Nevertheless, we can characterize ideal flow clutters by the absence of weakly MNI strong minors.

\begin{lem}
\label{lem:flow-clutter-minors}
For any flow clutter $\F$ the following are equivalent:
\begin{enumerate}[(i)]
\item $\F$ is ideal
\item $\F$ has no weakly MNI strong minor
\end{enumerate}
\end{lem}

The goal of this paper is to determine signed graphs that correspond to weakly MNI flow clutters, as their absence as strong minors characterizes idealness.
We first establish some simple properties.

\begin{lem}
\label{lem:parallel-edges}
Let $\F$ be weakly MNI. Then $G$ has no parallel edges.
\end{lem}

It is straightforward to see that flow clutters and minors of flow clutters obtained by contraction of $E^-$ do not correspond to degenerate projective planes.

\begin{lem} \label{lem:degenerate-proj-planes}
There is no degenerate projective plane that corresponds to a flow clutter. Similarly, no degenerate projective plane of order $k \geq 3$ corresponds to a clutter that is obtained from a flow clutter by contraction of all negative edges (purely as members of the clutter).
\end{lem}

Weakly MNI flow clutters have associated fractional vertices of the flow covering polyhedron.
In contrast to MNI clutters, however, some entries corresponding to negative edges may be $0$ and thus non-fractional.

\begin{lem} \label{lem:weakly-mni-frac-entries}
Let $\F$ be weakly MNI and suppose $x$ is a fractional vertex of $\mathsf P_{A(\F)}$. Then $0 < x_e < 1$ for all $e \in E^+$ and $x_e < 1$ for all $e \in E^-$.
\end{lem}

For any non-ideal flow clutter $\F$ and any fractional vertex $x$ of $\mathsf P_{A(\F)}$, let
\begin{align*}
E^-_0(x) = \{ e \in E^- \mid x_e = 0\}
\end{align*}
denote the subset of negative edges where $x$ takes the value $0$.
For weakly MNI flow clutters, contraction of the set $E^-_0(x)$ in the clutter yields an MNI minor.

\begin{lem} \label{lem:MNI-minor}
Let $\F$ be weakly MNI.
Then there exists a fractional vertex $x$ of $\mathsf P_{A(\F)}$ such that the clutter $\F / E^-_0(x)$ obtained by contraction of $E^-_0(x)$ is MNI.
In particular, if $x$ is a fractional vertex of $\mathsf P_{A(\F)}$ and $E^-_0(x) = E^-$, then $\F / E^-$ is MNI.
\end{lem}



\section{Idealness of Flow Clutters}
\label{sec:idealness}

In this section we discuss the idealness of flow clutters. We first employ the stronger notion of balancedness to characterize an instructive special case.

A square 0--1-matrix $A$ is called \emph{2-circulant} if there exist appropriate permutations of the rows and columns such that the permuted matrix $\tilde A$ is of the form
\begin{align*}
\tilde A =
\begin{pmatrix}
1 & 1 & & \\
& \ddots & \ddots & \\
& & 1 & 1 \\
1 & & & 1 \\
\end{pmatrix}.
\end{align*}
A 0--1-matrix $A$ is called \emph{balanced} if it does not contain any 2-circulant submatrix of odd size. Balanced matrices were introduced by \cite{Berge1972}, who also proved the following result.

\begin{thm}[\cite{Berge1972}] \label{thm:berge}
Any balanced 0--1-matrix is ideal.
\end{thm}

\subsection{Positive Trees}

We first consider the case that $G^+ = (V,E^+)$ is a tree.
Note that in this case problem \eqref{eq:correlation-clustering} remains NP-hard \citep{Garg1997}.

\begin{figure}
\center
(a)
\begin{tikzpicture}[scale=0.9,yshift=6ex]
\node (0) [style=vertex] at (0, 0) {};
\node (1) [style=vertex] at (-0.87, -0.5) {};
\node (2) [style=vertex] at (0, 1) {};
\node (3) [style=vertex] at (0.87, -0.5) {};

\draw (0) edge[thick] (1);
\draw (0) edge[thick] (2);
\draw (0) edge[thick] (3);
\draw (1) edge[dashed,bend left] (2);
\draw (1) edge[dashed,bend right] (3);
\draw (2) edge[dashed,bend left] (3);
\end{tikzpicture}
\hfill
(b)
\begin{tikzpicture}[scale=1,yshift=6ex]
\node (0) [style=vertex] at (0, 1) {};
\node (4) [style=vertex] at (-0.95, 0.31) {};
\node (1) [style=vertex] at (0.95, 0.31) {};
\node (3) [style=vertex] at (-0.59, -0.81) {};
\node (2) [style=vertex] at (0.59, -0.81) {};

\draw (0) edge[thick] (1);
\draw (1) edge[thick] (2);
\draw (2) edge[thick] (3);
\draw (3) edge[thick] (4);
\draw (4) edge[thick] (0);
\draw (0) edge[dashed] (2);
\draw (1) edge[dashed] (3);
\draw (2) edge[dashed] (4);
\draw (3) edge[dashed] (0);
\draw (4) edge[dashed] (1);
\end{tikzpicture}
\hfill
(c)
\begin{tikzpicture}[scale=0.9,yshift=6ex]
\node (0) [style=vertex] at (0.87, -0.5) {};
\node (1) [style=vertex] at (-0.87, -0.5) {};
\node (2) [style=vertex] at (0, 1) {};

\draw (0) edge[thick] (1);
\draw (1) edge[thick] (2);
\draw (2) edge[thick] (0);
\draw (1) edge[dashed,bend left] (2);
\draw (2) edge[dashed,bend left] (0);
\draw (0) edge[dashed,bend left] (1);
\end{tikzpicture}
\hfill
(d)
\begin{tikzpicture}[scale=0.8,yshift=6ex]
\node (0) [style=vertex] at (-0.1, 0) {};
\node (1) [style=vertex] at (1.1, 0) {};
\node (2) [style=vertex] at (-0.6, -0.8) {};
\node (3) [style=vertex] at (1.6, -0.8) {};
\node (4) [style=vertex] at (-0.5, 0.8) {};
\node (5) [style=vertex] at (1.5, 0.8) {};

\draw (0) edge[thick] (1);
\draw (1) edge[thick] (2);
\draw (2) edge[thick] (3);
\draw (0) edge[thick] (3);
\draw (0) edge[thick] (4);
\draw (1) edge[thick] (5);
\draw (2) edge[thick] (4);
\draw (3) edge[thick] (5);
\draw (0) edge[dashed] (2);
\draw (1) edge[dashed] (3);
\draw (4) edge[dashed] node[above] {$f$} (5);
\end{tikzpicture}
\caption[Forbidden strong minors]{Dashed lines depict negative edges, solid lines depict positive edges. (a) An odd flow-star with three positive edges. (b) The smallest non-ideal odd flow-circuit. (c) The only ideal odd flow-circuit. (d) Flow-split-$K_5$ with fractional vertex $x$ such that $x_f = 0$ and $x_e = 1/3$ for all $e \neq f$.}
\label{fig:non-ideal-minors}
\end{figure}
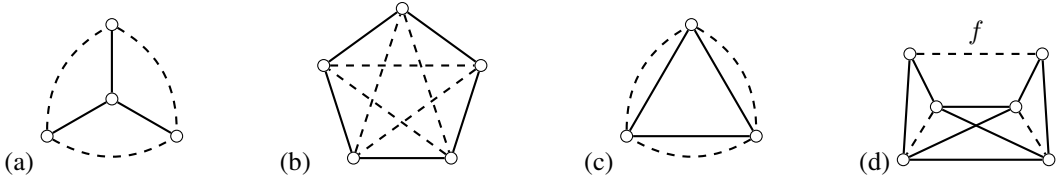

A \emph{flow-star} is a signed graph $S = (V_S, E_S^+, E_S^-)$ on $V_S = \{v_0, v_1, \dotsc, v_k\}$, where $k \geq 3$, such that $E_S^+ = \{ v_0 v_i \mid 1 \leq i \leq k\}$ and $E_S^- = \{ v_i v_{i+1} \mid 1 \leq i \leq k-1 \} \cup \{v_k v_1\}$. A flow-star is called odd if $k$ is odd. See Figure \ref{fig:non-ideal-minors} (a) for an example. Any odd flow-star $S$ defines a non-ideal flow clutter as the vector $x$ defined by $x_e = 1/2$ for $e \in E_S^+$ and $x_e = 0$ for $e \in E_S^-$ is a fractional vertex of the associated polyhedron.
Odd flow-star strong minors prohibit balancedness of the constraint matrix.

\begin{lem} \label{lem:positive-trees}
Let $G^+ = (V,E^+)$ be a tree. If $G$ has no odd flow-star strong minor, then $A(\F)$ is balanced.
\end{lem}

With Lemma~\ref{lem:positive-trees} we obtain the following result.

\begin{cor}
\label{cor:positive-trees}
Let $G^+ = (V,E^+)$ be a tree. Then $\F$ is ideal if and only if $G$ has no odd flow-star strong minor.
\end{cor}

\subsection{Positive Circuits}

Now we consider the case that $G^+=(V,E^+)$ is a circuit.
For weakly MNI flow clutters on positive circuits, we can show that the minor obtained from contraction of all negative elements is MNI.

\begin{lem} \label{lem:circuit-frac-vertex}
Let $G^+$ be a circuit and suppose that $\F$ is weakly MNI. There exists a fractional vertex $x$ of $\mathsf P_{A(\F)}$ such that $E^-_0(x) = E^-$.
\end{lem}

A \emph{flow-circuit} is a signed graph $C = (V_C,E_C^+,E_C^-)$ on $V_C = \{v_1, \dotsc, v_k\}$, where $k \geq 3$, such that $E_C^+ = \{ v_i v_{i+1} \mid 1 \leq i \leq k-1 \} \cup \{v_k v_1\}$ and $E_C^- = \{ v_i v_{i+2} \mid 1 \leq i \leq k-2 \} \cup \{v_{k-1} v_1, v_k v_1\}$.
A flow-circuit is called odd if $k$ is odd. Any odd flow-circuit $C$ with $\abs{E_C^+} \geq 5$ defines a non-ideal flow clutter as the vector $x$ defined by $x_e = 1/2$ for $e \in E_C^+$ and $x_e = 0$ for $e \in E_C^-$ is a fractional vertex of the associated polyhedron.
See Figure \ref{fig:non-ideal-minors} (b) for an example.
By another application of Lehman's Theorem we can characterize idealness of $\F$ by forbidden odd flow-circuit strong minors as follows.

\begin{thm} \label{thm:positive-circuit}
Let $G^+$ be a circuit. Then the clutter $\F$ is ideal if and only if $G$ has no odd flow-circuit strong minor $C = (V_C,E_C^+,E_C^-)$ with $\abs{E_C^+} \geq 5$.
\end{thm}

\subsection{General Graphs}

In the case of positive trees and positive circuits we exploited the fact that for weakly MNI flow clutters $\F$, any fractional vertex $x$ of $\mathsf P_{A(\F)}$ such that $\F / E^-_0(x)$ is MNI satisfies $E^-_0(x) = E^-$.
In this section we show for general graphs, if $E^-_0(x) \neq E^-$, then the core of the corresponding MNI clutter satisfies a property that is encountered with very few known cores.

The core $\overline \C$ of an MNI clutter $\C$ is called \emph{fat} if $cb - n + 1 \geq 3$, where $c,b$ and $n$ are the constants from Theorem \ref{thm:lehman}. Currently, only three distinct fat cores of MNI clutters are known \citep{Cornuejols2009}, the clutter $F_7$ of lines of the Fano plane (which is its own blocker), the clutter $\tau(K_5)$ of triangles of $K_5$ and its blocker. In each case it holds that $cb - n + 1 = 3$.

%

\begin{lem}
\label{lem:fat-core}
Let $\F$ be a weakly MNI flow clutter.
Let $x$ be a fractional vertex of $\mathsf P_{A(\F)}$ such that $\C = \F / E^-_0(x)$ is MNI and $E^-_0(x) \neq E^-$.
Then the core $\overline \C$ is fat.
\end{lem}

The instance depicted in Figure~\ref{fig:non-ideal-minors} (d), which we call \emph{flow-split-$K_5$}, is weakly MNI.
Let~$\F$ be the corresponding flow clutter.
It has an associated fractional vertex $x$ such that $x_f = 0$ and $x_e = 1/3$ for all $e \neq f$.
Further, the core of the MNI minor $\F / f$ is isomorphic to~$\tau(K_5)$.
To see this, take $K_5$, split an arbitrary vertex into two vertices with two neighors each, connect them by $f$ and sign the resulting graph appropriately.
The other known fat cores of MNI clutters do not arise from flow clutters, as we point out below.
This shows that any weakly MNI flow clutter that arises from a vertex $x$ of $\mathsf P_{A(\F)}$ with $E^-_0(x) \neq E^-$ and is different from the flow-split-$K_5$ would imply the existence of an unknown MNI clutter with a fat core.

\begin{lem}
\label{lem:nonisomorphic-cores}
Let $\F$ be a weakly MNI flow clutter and $x$ be a fractional vertex of $\mathsf P_{A(\F)}$ such that the clutter $\C = \F / E^-_0(x)$ is MNI.
If $E^-_0(x) \neq E^-$, then~$\C$ is neither isomorphic to the blocker of $\tau(K_5)$ nor to the clutter~$F_7$.
\end{lem}

\subsection{Discussion}

In the proof for positive circuits, we consider the MNI minor of the flow clutter that consists of the corresponding positive paths.
In fact, for both positive trees and positive paths we exploit the one-to-one correspondence between flows and its positive subpaths.
Therefore it may seem more useful to consider the clutter of positive paths to start with and instead characterize the MNI minors of that clutter.
After all, that clutter corresponds to the minimum multicut formulation, where one seeks to cover a set of paths between pairs of terminal vertices.
In this section, we discuss this issue by presenting two arguments why our approach is favorable nonetheless.

i. There is a one-to-one correspondence between strong minors of flow clutters and strong minors of signed graphs.
If instead we consider the clutter of terminal paths in an unsigned graph, then the minors of that clutter are not in one-to-one correspondence with the minors of the graph, due to the fact that the graph minors do not carry any information about the terminal pairs.
In the language of minimum multicut, one needs to consider the union of the supply and demand graph.

\begin{figure}[!h]
\center
\begin{tikzpicture}[scale=1]
\node (0) [vertex] at (0, 0) {};
\node (1) [vertex] at (-0.71, 0.71) {};
\node (2) [vertex] at (-0.71, 1.71) {};
\node (3) [vertex] at (0, 2.42) {};
\node (4) [vertex] at (1, 2.42) {};
\node (5) [vertex] at (1.71, 1.71) {};
\node (6) [vertex] at (1.71, 0.71) {};
\node (7) [vertex] at (1, 0) {};

\draw (0) edge (1);
\draw (1) edge (2);
\draw (2) edge (3);
\draw (3) edge (4);
\draw (4) edge (5);
\draw (5) edge (6);
\draw (6) edge (7);
\draw (7) edge (0);
\draw (0) edge[dashed] (3);
\draw (1) edge[dashed] (4);
\draw (2) edge[dashed] (5);
\draw (3) edge[dashed] (6);
\draw (4) edge[dashed] (7);
\draw (5) edge[dashed] (0);
\draw (6) edge[dashed] (1);
\draw (7) edge[dashed] (2);
\end{tikzpicture}
\caption[An advantage of weakly MNI strong minors]{Dashed lines depict negative edges, solid lines depict positive edges. The depicted instance is not weakly MNI, but has an odd flow-circuit strong minor. If instead the negative edges are interpreted as terminal pairs, then the corresponding clutter of terminal paths is already MNI.}
\label{fig:circulant-instance}
\end{figure}
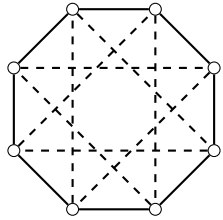

ii. There are fewer minors that need to be forbidden.
This can be seen by considering instances defined on positive circuits.
Note that the odd flow-stars and odd flow-circuits are instances of so-called \emph{circulant} clutters, which can be represented by circulant 0--1-matrices.
In fact, both structures correspond to 2-circulant matrices.
\cite{Cornuejols1994} provide a complete list of all other MNI circulant clutters.
If we consider the clutter of terminal paths defined on a circuit, then it turns out that it does not suffice to only forbid minors that correspond to 2-circulant clutters.
For example, consider the instance depicted in Figure~\ref{fig:circulant-instance}.
It is not weakly MNI, but non-ideal as it has an odd flow-circuit strong minor.
If instead the negative edges are regarded as terminal pairs, then the core of the corresponding clutter of terminal paths is the MNI circulant clutter~$\C^3_8$, which consists of 8 members with 3 elements each.

\section{Conclusion and Future Work}

We introduced the class of flow-partitionable signed graphs, which allow to solve the associated correlation clustering problem efficiently via its cycle relaxation.
Furthermore, we proposed to characterize flow-partitionable signed graphs in terms of forbidden minors.
The forbidden minors we found include two infinite classes that provide a complete characterization for the cases when the positive subgraph is either a circuit or a tree.
For the general case, we provide another forbidden minor.
Our results indicate that a complete characterization may be difficult.
Therefore, the restriction to other special cases seems reasonable in future work.
In particular, we think the study of planar signed graphs is promising for three reasons.
(i) The forbidden flow-circuit minors and the flow-split-$K_5$ minor that we presented are non-planar.
(ii) Planar graphs are only moderately less sparse than series-parallel graphs (no $K_5$ and no $K_{3,3}$ minor opposed to no $K_4$ minor) and series-parallel signed graphs are flow-partitionable \citep{Chopra1994}.
(iii) Planar correlation clustering instances admit a PTAS \citep{Klein2015} whereas for general graphs any constant factor approximation is NP-hard under the Unique Games Conjecture \citep{Chawla2006}.
We thus raise the question: Is a signed graph $G$ flow-partitionable if $G$ is planar and has no odd flow-star strong minor?

\bibliographystyle{abbrvnat}
\bibliography{main}


\appendix

\section{Proofs}
\label{sec:proofs}

\subsection{Section \ref{sec:ideal-clutters}}

\paragraph{Lemma \ref{lem:flow-clutter-minors}}

\begin{proof} 
Every weakly MNI strong minor has an MNI minor, so if $\F$ is ideal, then it cannot have a weakly MNI strong minor.
Conversely, suppose $\F$ has no weakly MNI strong minor, but some MNI minor $\F'$.
Then, since the contraction and deletion operations are commutative, this minor
is obtained from a strong minor $\F''$ by a series of contractions of negative edges.
Thus, $\F''$ is weakly MNI, which is a contradiction.
\end{proof}

\paragraph{Lemma \ref{lem:parallel-edges}}

\begin{proof}
Let $f,g \in E$ be a pair of parallel edges.
First assume that $f, g \in E^-$ are both negative.
Let $x$ be any vertex of $\mathsf P_{A(\F)}$.
If $x_f < x_g$, then $x = 1/2 (y + z)$, where $y, z \in \mathsf P_{A(\F)}$ agree with $x$ except for $y_g = x_f$ and $z_g = 2x_g - x_f$.
Thus, it must hold that $x_f = x_g$.
Now suppose $x$ is a fractional vertex of $\mathsf P_{A(\F)}$.
Let $x_{\backslash f} \in \mathsf P_{A(\F\backslash f)}$ denote the vector obtained from $x$ by setting $x_f = 1$.
Since $\F$ is weakly MNI, $x_{\backslash f}$ can be written as a convex combination of 0-1-vectors $y^k \in \mathsf P_{A(\F\backslash f)}$.
Every $y^k$ is a vertex, so $y^k_g = y^k_f = 1$. This implies that $x_f = x_g = 1$, which is a contradiction.

If $f,g \in E^+$ are both positive, then the proof is completely analogous.
Otherwise (w.l.o.g.) $f \in E^+$ and $g \in E^-$, so the cycle induced by $f$ and $g$ is a flow.
Let $x$ be any vertex of $\mathsf P_{A(\F)}$. We have $x_f + x_g \geq 1$. Assume that $x_f + x_g > 1$. Since $x$ is a vertex, neither $x_f$ nor $x_g$ can be decreased without violating an inequality of $\mathsf P_{A(\F)}$.
Therefore, there exists a flow $C$ containing $f$ and a path $P$ that induces a flow with $g$ such that $x_P + x_g = 1$ and $x_{C \setminus \{f\}} + x_f = 1$. Further, the set $P \cup C \setminus \{f\} $ contains a flow.
This implies that
\begin{align}
& 1 \leq x_{C \setminus \{f\}} + x_P = 2 - x_f - x_g \\ \implies & x_f + x_g \leq 1,
\end{align}
which is a contradiction.
Hence, it holds that $x_f + x_g = 1$.
Now suppose $x$ is a fractional vertex of $P_{A(\F)}$. Since $\F$ is weakly MNI, the vector $x_{\backslash f}$ can be written as a convex combination of 0-1-vectors $y^k \in \mathsf P_{A(\F \backslash f)}$.
Every $y^k$ is a vertex, so they satisfy $y^k_f = 1$ and $y^k_g = 0$.
This implies that also $x_g = 0$ and $x_f = 1$, which is a contradiction.
\end{proof}

\paragraph{Lemma \ref{lem:degenerate-proj-planes}}

\begin{proof}
Let $\D = \big \{ \{1, \dotsc, k\}, \{0,1\}, \dotsc, \{0,k\} \big \}$ be a degenerate projective plane of order $k$. If $0 \in E^-$, then $1, \dotsc, k \in E^+$, which implies that the first member of $\D$ has no negative edge. If $0 \notin E^+$, then $1, \dotsc, k \in E^-$, which implies that the first member of $\D$ has no positive edge.

Now, suppose $\D$ consists of the edge sets of the positive paths of a flow clutter. Since $\{1, \dotsc, k\}$ is a path (and not a cycle) of length $k \geq 3$ one of the members $\{0,1\}, \dotsc, \{0,k\}$ cannot represent a connected path.
\end{proof}

\paragraph{Lemma \ref{lem:weakly-mni-frac-entries}}

\begin{proof}
Every vertex of $\mathsf P_{A(\F)}$ satisfies $0 \leq x_e \leq 1$ for all $e \in E$. If $x_e = 1$ for any $e \in E$, then delete $e$ to obtain a strong minor of $\F$ that is non-ideal. If $x_e = 0$ for $e \in E^+$, then contract $e$ to obtain a strong minor of $\F$ that is non-ideal.
\end{proof}

\paragraph{Lemma \ref{lem:MNI-minor}}

\begin{proof}
Take some fractional vertex $x$ of $\mathsf P_{A(\F)}$ and suppose the clutter $\C = \F / E^-_0(x)$ is not MNI.
If there exists an element $e \in E^+ \cup E^- \setminus E^-_0(x)$ such that $\C \backslash e$ is non-ideal, then, since $\C \backslash e = (\F \backslash e) / E^-_0(x)$, it follows that $\F \backslash e$ is non-ideal, so $\F$ is not weakly MNI.
A similar argument shows that there cannot be any $e \in E^+$ such that $\C / e$ is non-ideal.
Thus, there exists some $f \in E^- \setminus E^-_0(x)$ such that $\C / f$ is non-ideal. Hence, $P_{A(\F)}$ has a fractional vertex $y$ with $y_e = 0$ for all $e \in E^-_0(x)$ and $y_f = 0$. Replacing $x$ by $y$ and repeating this argument eventually yields the desired $x$, as $E^-$ is finite.
The argument further shows that if $E^-_0(x) = E^-$, then $\F / E^-$ must already be MNI.
\end{proof}

\subsection{Section \ref{sec:idealness}}

\paragraph{Lemma \ref{lem:positive-trees}}

\begin{proof}
For contraposition, assume $A(\F)$ is not balanced, so it has a 2-circulant submatrix $B$ of odd order. Note that, since $G^+$ is a tree, every negative edge induces exactly one flow. Thus, the columns of $B$ correspond to positive edges only.
We construct an odd flow-star strong minor of $G$. First, delete all edges that do not correspond to any flow associated with the rows of $B$. Then, contract all positive edges from $G$ that do not correspond to any column of $B$.
This yields a minor $S = (V_S, E_S^+, E_S^-)$ of $G$. 
Clearly, since $G^+$ is a tree, $S^+$ is also a tree. Further, we never contract any parallel edge which shows that $S$ is a strong minor.
The structure of $B$ implies that the graph $S$ has an odd number of flows of length three. The flows can be cyclically ordered such that every adjacent pair of flows shares a positive edge. Hence, it follows that $S$ must be an odd flow-star.
\end{proof}

\paragraph{Corollary \ref{cor:positive-trees}}

\begin{proof}
If $G$ has no odd flow-star strong minor, then $\F$ is ideal by Lemma \ref{lem:positive-trees} and Theorem~\ref{thm:berge}. Conversely, if $G$ has an odd flow-star strong minor, then $\F$ cannot be ideal, since any odd flow-star is non-ideal and idealness is preserved under taking minors.
\end{proof}

\paragraph{Lemma \ref{lem:circuit-frac-vertex}}

\begin{proof}
Let $x$ be a fractional vertex of $\mathsf P_{A(\F)}$ such that $\F / E^-_0(x)$ is MNI and assume that $0 < x_f < 1$ for some $f \in E^-$. 
If $\F / E^-_0(x)$ is a degenerate projective plane of order $k \geq 2$, then it has a member of size two that contains $f$. Thus, $f$ is parallel to some positive edge, which is a contradiction to Lemma~\ref{lem:parallel-edges}.
Hence, the $\F / E^-_0(x)$ is not a degenerate projective plane and, by Theorem~\ref{thm:lehman}, the flows that share $f$ have size two, since $f$ induces exactly two flows. This implies that $E^+$ has size two and therefore $\F$ is ideal, which is a contradiction.
\end{proof}

\paragraph{Theorem \ref{thm:positive-circuit}}

\begin{proof}
If $G$ has an odd flow-circuit strong minor $C$ with $\abs{E_C^+} \geq 5$, then $\F$ cannot be ideal as idealness is preserved under taking minors.

Conversely, suppose $\F$ is non-ideal and weakly MNI. Consider the minor $\P = \F / E^-$ obtained by contraction of all negative elements in $\F$. Since $\F$ is weakly MNI, the clutter~$\P$ is MNI (Lemmas~\ref{lem:MNI-minor} and \ref{lem:circuit-frac-vertex}) and consists of the edge sets of paths in $G^+$ associated with the flows in $\F$. As $\P$ is not a degenerate projective plane (Lemma \ref{lem:degenerate-proj-planes}), Theorem \ref{thm:lehman} implies that there are $n = \abs{E^+}$ minimum members of $\overline \P$ that all have some constant length $p$. Let~$\M$ be the blocker of $\P$ so the $n$ members of $\overline \M$ have some constant size $m$.
By Theorem~\ref{thm:lehman}, it holds that $pm \geq n + 1$. Further, we have $2 \leq p \leq \frac{n}{2}$ as $G^+$ is a circuit and the members of $\overline \P$ have minimum length. 
As the $n$ paths in $\overline \P$ of length $p$ can be arranged cyclically, it is apparent that $m = \big \lceil \frac{n}{p} \big \rceil$. Now, if $p = \frac{n}{2}$, then $m = 2$ and thus $pm = n < n + 1$. Therefore, we must have that $p \leq \frac{n-1}{2}$. In particular, every negative edge $f \in E^-$ corresponds to exactly one path $P_f \in \overline \P$ and $\overline \P = \P$.

Assume that $p \geq 3$. We show that this leads to a contradiction and thus $p = 2$.
Take some $M \in \overline \M$. By Theorem \ref{thm:lehman}, there exists a unique $P_f \in \P$ for $f \in E^-$ such that $\abs{P_f \cap M} = pm - n + 1$ and $\abs{P \cap M} = 1$ for all $P \in \P$ with $P \neq P_f$.
We define a vector $x \in \R^{E^+ \cup E^-}$ by
\begin{align}
x_e = \begin{cases} \frac{1}{p-1} & e \in E^+ \setminus M \\ 0 & e \in M \\ \frac{pm-n}{p-1} & e = f \\ 0 & e \in E^- \setminus f. \end{cases}
\end{align}
The vector $x \geq 0$ is constructed such that for every $P \in \P$, the corresponding flow covering inequality of $\mathsf P_{A(\F)}$ is tight. Indeed, it holds that
\begin{align}
x(P_f \cup \{f\}) = \frac{p-pm+n-1}{p-1} + \frac{pm-n}{p-1} = \frac{p-1}{p-1} = 1
\end{align}
and, for all $e \in E^-$ with $e \neq f$ that
\begin{align}
x(P_e \cup \{e\}) = \frac{p-1}{p-1} + 0 = 1.
\end{align}
Feasibility of $x$ is clear, since each path in $E^+$ that corresponds to a flow but is not of minimum size contains at least two distinct members of $\P$.
Therefore, it holds that $x$ is a vertex of $\mathsf P_{A(\F)}$. Further, since $p \geq 3$, the vertex $x$ is fractional, which implies that $\F$ cannot be weakly MNI. Thus, it must hold that $p = 2$.

This shows that, as $\P$ is MNI, the matrix $A(\P)$ must be an odd 2-circulant matrix. Moreover, it holds that $\abs{E^+} \geq 5$ as otherwise $\F$ would be ideal. Hence, the signed graph $G$ is an odd flow-circuit with $\abs{E^+} \geq 5$.
\end{proof}

\paragraph{Lemma \ref{lem:fat-core}}

\begin{proof}
Since $E^-_0(x) \neq E^-$ there exists some $f \in E^-$ such that $0 < x_f < 1$.
Assume that $\C$ is a degenerate projective plane of order $k \geq 2$.
Then $f$ is parallel to some positive edge, which is a contradiction to $\F$ being weakly MNI.

Thus, the clutter $\C$ is not a degenerate projective plane. Let $\overline \C = \{C_i\}_i$ and $\overline \B = \{B_i\}_i$ denote the cores of $\C$ and its blocker $\B$.
By Theorem \ref{thm:lehman}, there are $i, j \in [n]$ with $i \neq j$ such that $f \in C_i \cap B_i$ and $f \in C_j \cap B_j$.
Further, it holds that $C_i \cap C_j = \{f\}$.
Therefore, the set of edges $(C_i \cup C_j) \setminus \{f\}$ induces a positive cycle.
Since $\abs{C_i \cap B_i} = cb - n + 1 \geq 2$, there exists some $g \in E^+$ such that $g \in C_i \cap B_i$.
There is another flow $C_k \in \overline \C$, $k \neq i,j$ such that $C_i \cap C_k = \{g\}$.
Now, since $(C_i \cup C_j) \setminus \{f\}$ is a positive cycle, the set $(C_k \cup C_i \cup C_j) \setminus \{f,g\}$ contains a flow, which must be covered by $B_i$.
Hence, there exists another edge $h \in (C_i \cup C_j) \setminus \{f,g\} \cap B_i$. It must hold that $h \in C_i$ as $C_j \cap B_i = \{f\}$.
This shows that $\{f,g,h\} \subseteq C_i \cap B_i$ and thus $cb - n + 1 \geq 3$.
\end{proof}

\paragraph{Lemma \ref{lem:nonisomorphic-cores}}

\begin{proof}
First, assume that $\C$ is isomorphic to the blocker of $\tau(K_5)$.
Then $\overline \C$ has 10 members with~4 elements each.
By assumption, there is some $f \in E^-$ that is contained in 4 members of $\C$.
Since all members are of minimum size and no edges may be parallel, at least 8 positive edges are needed to form the flows that include $f$ and they are arranged as in Figure~\ref{fig:non-ideal-minors}~(d) (without the other two negative edges).
Furthermore, it is clear that $\overline \C$ cannot contain this substructure, as its members are composed of only 10 distinct elements.

Next, assume $\C$ is isomorphic to $F_7$ and take some $f \in E^- \setminus E^-_0(x)$. As $x_f > 0$, one element of $F_7$ corresponds to $f$.
Thus, since for any other element $e$ of $F_7$ there is a member of $F_7$ that contains both $e$ and $f$, it holds that~$e \in E^+$.
It follows that $G^+$ consists of three edge-disjoint paths of length two connecting the endpoints of $f$.
Moreover, there is another member of~$F_7$ that contains one edge from each such path and forms a path itself.
This is impossible and thus the edge $f$ cannot exist.
\end{proof}

\end{document}